\documentclass[twoside,reqno,A4]{amsart}

\numberwithin{equation}{section}

\newcommand{\qdn}{\hspace*{-1.5mm}}
\newcommand{\qqdn}{\hspace*{-2.5mm}}

\newcommand{\+}{&\qqdn}%

\newcommand{\mb}[1]{\mathbb{#1}}

\newcommand{\ffnkK}[4]{\left[\qdn\ba{#1}#3\\[2mm]#4\ea{\!\:#2}\right]}

\newcommand{\ffnk}[4]{\left[\qdn\ba{#1}#3\\[2mm]#4\ea{\!;\:#2}\right]}

\newcommand{\nnm}{\nonumber}
\newcommand{\be}{\begin{equation}}
\newcommand{\ee}{\end{equation}}
\newcommand{\ba}{\begin{array}}
\newcommand{\ea}{\end{array}}
\newcommand{\bmn}{\begin{eqnarray}}
\newcommand{\emn}{\end{eqnarray}}
\newcommand{\bnm}{\begin{eqnarray*}}
\newcommand{\enm}{\end{eqnarray*}}
\newcommand{\bln}{\begin{subequations}}
\newcommand{\eln}{\end{subequations}}

\newtheorem{thm}{Theorem}

\newcommand{\bbtm}[4]{\bibitem{kn:#1}{#2,}~\emph{#3,}~{#4.}}
\newcommand{\cito}[1]{\cite{kn:#1}}
\newcommand{\citu}[2]{\cite[#2]{kn:#1}}

\begin{document} 
\title{Comments on "New generating relations for products of two Laguerre polynomials"}
\author{Xiaoxia Wang $^{a, *}$, Arjun K. Rathie $^b$}
\dedicatory{$^a$ Department of Mathematics, Shanghai University, Shanghai, 200444, P. R. China;\\[1mm]
$^b$ Department of Mathematics, School of Mathematical and Physical Sciences, Central University of Kerala, \\[1mm]
Riverside Transit Campus, Padennakkad P. O. Nileshwar Kasaragod-671 328, Kerala State, India.}

\thanks{$^{*}$ Corresponding author. \\
E-mail addresses: xiaoxiawang@shu.edu.cn (X. Wang), akrathie@rediffmail.com (A. K. Rathie).}



\maketitle\thispagestyle{empty}
\markboth{X. Wang and A.K. Rathie}
{Comments on "New generating relations for products of two Laguerre polynomials"}
\begin{center}\parbox{100mm}{\footnotesize
\emph{Abstract.} By utilizing a two-dimensional extension of a very general series transform given by Bailey,
Exton [Indian J. pure appl. Math. 24 (6) (1993), 401-408] deduced a very general double
generating relation of a product of a pair of Laguerre polynomials
and obtained a number of useful relations with elementary functions,
Bessel functions, Hermite polynomials and single series expansions
of pairs of Laguerre polynomials. Unfortunately, some of the results
given by Exon contain errors and thus this is the aim of this short note
to provide the corrected form of these results.\\

\emph{2000 AMS Subject Classification}: 33C20; 33C05; 33B20.\\

\emph{Keywords}: Generating relations;
          Laguerre polynomial;
          Bessel function;
          Hermite polynomial.}

\end{center}

\section{\bf{Introduction}}

We recall with the definition of generalized hypergeometric function \cito{rainville} with $p$ numerator and $q$ denominator parameters by
\bmn\label{def}
{_p F_q}
\ffnkK{ccc}{z}{\alpha_1,\+\cdots,\+\alpha_p;}
              {\beta_1,\+\cdots,\+\beta_q;}
={_p F_q}[(\alpha); (\beta); z]
=\sum_{n=0}^{\infty}
\frac{(\alpha_1,n)\cdots(\alpha_p,n)}
     {(\beta_1,n)\cdots(\beta_q,n)}
 \frac{z^n}{n!}
 =\sum_{n=0}^{\infty}
\frac{((\alpha),n)}{((\beta),n)}\frac{z^n}{n!},
\emn
where $(\alpha, n)$ denotes the well known Pochhammer symbol
(or the shifted factorial, since $(1,n)=n!$) defined for complex number $\alpha$ by
\bmn
(\alpha, n)
=\begin{cases}
\alpha(\alpha+1)\cdots(\alpha+n-1),\quad \+n\in\mb{N};\\[2mm]
1,\+n=0.
\end{cases}
\emn
Using the fundamental function relation $\Gamma(\alpha+1)=\alpha\Gamma(\alpha)$, $(\alpha, n)$ can be written in the form
\bmn\label{gamma}
(\alpha, n)=\frac{\Gamma(\alpha+n)}{\Gamma(\alpha)}, \quad \quad \quad (n\in\mb{N}\cup \{0\}),
\emn
where $\Gamma$ is the well known Gamma function. For more detail about
convergence of this function, we refer to \cito{rainville}.

On the other hand, we recall the double hypergeometric function which is
defined and introduced by Kamp\'{e} de F\'{e}riet and subsequently
abbreviated by Bunchnall and Chaundy \cito{sri-kar}.
Here, we present a slightly modified notation given by
Srivastava and Panda \citu{sri-pan}{p.423, Eq.(26)} as follows.

\bmn \label{def-de}
{{F}_{g: \:c; \:d}^{h: \:a; \:b}}
\ffnkK
{ccc}{x,y}{(H_h):\+(A_a); \+(B_b);}{(G_g):\+(C_c); \+(D_d);}
=\sum_{m,n\geq0}
\frac{((H_h),m+n)((A_a),m)((B_b),n)}{((G_g),m+n)((C_c),m)((D_d),n)}
\frac{x^m}{m!}\frac{y^n}{n!}.
\emn
For more detail about the convergence of this function, we refer to \cito{sri-pan}.

The Laguerre polynomials have been researched in various branches
 of pure and applied mathematics \cite{kn:Erdelyi-2, kn:slater}, which can be expressed by the confluent hypergeometric function as
\bmn\label{def-L}
L_n^{(a)}(x)=\frac{(a+1,n)}{n!}{_1F_1}\ffnk{c}{x}{-n}{a+1}.
\emn

In 1974, Exton \cito{exton-2} obtained the well known Bailey's transform in two dimension
in the following theorem.
\begin{thm}
If
\bmn\label{t-1}
\beta_{m, n}=\sum_{p=0}^m\sum_{q=0}^n\alpha_{p, q}\mu_{m-p, n-q}\nu_{m+p, n+q}
\emn
and
\bmn\label{t-2}
\gamma_{m, n}=\sum_{p=m}^\infty\sum_{q=n}^\infty\delta_{p, q}\mu_{p-m, q-n}\nu_{p+m, q+n},
\emn
then
\bmn\label{t-3}
\sum_{m, n=0}^\infty\alpha_{m, n}\gamma_{m, n}=\sum_{m, n=0}^\infty\beta_{m, n}\delta_{m, n}.
\emn
\end{thm}
Here, it is understood that $\alpha, \delta, \mu$ and $\nu$ are functions of $p$ and $q$
only and all the series concerned are either convergent or terminating.

By utilizing the above theorem and \eqref{t-1}, \eqref{t-2} and \eqref{t-3} are appropriately
specified, Exton \cito{exton} obtained the following very general double generating function

\bmn \label{main}
\+\+\sum_{m,n=0}^\infty\frac{((d),m+n)x^ms^n}{((g),m+n)m!n!}\nnm
{_1F_1}\ffnk{c}{-y}{-m}{p}{_1F_1}\ffnk{c}{-t}{-m}{p'}\\[2mm]
\+\+=\sum_{m,n=0}^\infty\frac{((d),m+n){(xy)}^m{(st)}^n}{((g),m+n)(p,m)(p',n) m!n!}
\:F_{G:0;0;}^{D:0;0;}\ffnkK{ccc}{x,s}{(d)+m+n:-;-;}{(g)+m+n:-;-;},
\emn
and further as a simple consequence of the binomial theorem, the inner double series on the
right-hand side of \eqref{main} is immediately reduce to a single series.
Also, if the confluent hypergeometric functions on the left-hand side of \eqref{main}
are replaced by their representations as Laguerre polynomials
(of course, changing $y$ to $-y$ and $t$ to $-t$ and using \eqref{def-L}),
we arrive at the following result.
\bmn \label{main-2}
\+\+\sum_{m,n=0}^\infty\frac{((d),m+n)x^ms^n}{((g),m+n)(p, m)(p',n)}\nnm
L_m^{p-1}(y)L_n^{p'-1}(t)\\[2mm]
\+\+=\sum_{m,n=0}^\infty\frac{((d),m+n){(-xy)}^m{(-st)}^n}{((g),m+n)(p,m)(p',n) m!n!}
\:{_D F_G}\ffnk{c}{x+s}{(d)+m+n}{(g)+m+n}.
\emn
This is a two-dimensional very general generating relation of a pair of Laguerre polynomials.

Exton in his paper \cito{exton} deduced a number of useful relations with
elementary functions, Bessel functions, Hermite polynomials and single series
expansions of pairs of Laguerre polynomials by utilizing \eqref{main-2}.
Unfortunately, some of the results given by Exton contain errors.

The remainder part of this short note is organized as follows. In section 2,
Exton's general result \eqref{main} will be established by another method.
In section 3, we will list Exton's results in corrected form.

\section{Another proof of \eqref{main}:}

In order to establish \eqref{main}, we proceed as follows. Denoting the
right-hand side of \eqref{main} by $\mathbf{S}$, expressing the
double hypergeometric function with the help of the definition \eqref{def-de}, we have
\bnm
\mathbf{S}\+=\+\sum_{m=0}^\infty\sum_{n=0}^\infty\sum_{u=0}^\infty\sum_{v=0}^\infty
\frac{((d),m+n)((d)+m+n,u+v)x^{m+u}s^{n+v} y^m t^n}
{((g),m+n)((g)+m+n,u+v)(p,m)(p',n) m!n! u!v!}\\
\+=\+\sum_{m=0}^\infty\sum_{n=0}^\infty\sum_{u=0}^\infty\sum_{v=0}^\infty
\frac{((d), m+n+u+v)x^{m+u}s^{n+v} y^m t^n}
{((g),m+n+u+v)(p,m)(p',n) m!n! u!v!},
\enm
where, we have applied the elementary relation $(a,m)(a+m,n)=(a, m+n)$.
Performing the transformations $u\to u-m$ and $v \to v-n$ on the above identity and then applying
the transformation $(m-n)!=(-1)^n m!/(-m, n)$, we get
\bnm
\mathbf{S}\+=\+\sum_{u=0}^\infty\sum_{v=0}^\infty\sum_{m=0}^u\sum_{n=0}^v
\frac{((d),u+v)x^{u}s^{v} y^m t^n}{((g), u+v)(p,m)(p',n) (u-m)!(v-n)!m!\:n!}\\
\+=\+\sum_{u=0}^\infty\sum_{v=0}^\infty\sum_{m=0}^u\sum_{n=0}^v
\frac{((d), u+v)(-u,m)(-v,n) x^{u}s^{v} (-y)^m (-t)^n}{((g), u+v)(p,m)(p',n) m!\:n!\: u!\:v!}\\
\+=\+\sum_{u=0}^\infty\sum_{v=0}^\infty\frac{((d), u+v)x^{u}s^{v}}{((g),u+v)u!\:v!}
\sum_{m=0}^u\sum_{n=0}^v\frac{(-u,m) (-v,n) (-y)^m (-t)^n}{(p,m)(p',n) m!\:n! }\\
\+=\+\sum_{u=0}^\infty\sum_{v=0}^\infty\frac{((d), u+v)x^{u}s^{v}}{((g), u+v)u!\:v!}
\sum_{m=0}^u\frac{(-u,m)}{(p,m) m! }(-y)^m
\sum_{n=0}^v\frac{(-v,n)}{(p',n) n! }(-t)^n\\
\+=\+\sum_{u=0}^\infty\sum_{v=0}^\infty\frac{((d), u+v)x^{u}s^{v}}{((g), u+v)u!v!}
{_1F_1}\ffnk{c}{-y}{-u}{p}{_1F_1}\ffnk{c}{-t}{-v}{p'}.
\enm
Let $u \to m $ and $v \to n$, the above identity arrives at the left-hand side of \eqref{main}.
This completes the proof of \eqref{main}.

\section{\bf{Exton's results in corrected form}}
In \cito{exton}, we find that some of the results are not correct. In this section,
we will present the corrected form of these results by the method applied by Exton
with specializing the parameters in the main transformation \eqref{main-2}.

Exton's result (3.3) should be read as
\bnm
\sum_{m,n=0}^\infty\frac{((d),m+n)(-1)^nx^{m+n}}{((g),m+n)(p,m)(p',n)}
L_{m}^{(p-1)}(y)L_{n}^{(p'-1)}(-y)\\
=_{D+2}F_{G+3}\ffnkK{ccc}{-4xy}{(d),\+(p+p'-1)/2,\+(p+p')/2;}{(g),\+p,\quad p',\+p+p'-1;};
\enm

Exton's result (3.8) should be read as
\bnm
\sum_{m,n=0}^\infty\frac{(p',m+n)(p+p'-1,m+n)(-1)^n x^{m+n}}
{((p+p'-1)/2,m+n)((p+p')/2,m+n)(p,m)(p',n)}\\
L_{m}^{(p-1)}(y)L_{n}^{(p'-1)}(-y)
=\Gamma(p)(2\sqrt{(xy)})^{1-p}J_{p-1}(4\sqrt{(xy)});
\enm

Exton's result (3.11) should be read as
\bnm
\sum_{m,n=0}^\infty\frac{(p,m+n)(p',m+n)(-1)^n x^{m+n}}
{(p+p',m+n)(p,m)(p',n)}L_{m}^{(p-1)}(-y)L_{n}^{(p'-1)}(y)\\
=\Gamma((p+p')/2)e^{2xy}(xy)^{1-p/2-p'/2}I_{p/2+p'/2-1}(2xy);
\enm

Exton's result (3.12) should be read as
\bnm
\sum_{m,n=0}^\infty\+\+\frac{(p,m+n)(p',m+n)(-1)^n x^{m+n}}
{(p,m)(p',n)}L_{m}^{(p-1)}(-y)L_{n}^{(p'-1)}(y)\\
\+=\+{_2F_1}((p+p'-1)/2,(p+p')/2;p+p'-1;4xy)\\
\+=\+(1-4xy)^{-1/2}[1/2+(1-4xy)^{1/2}/2]^{2-p-p'};
\enm

Exton's result (3.13) which is obtained by setting $p'=2-p$ in the above identity actually should be read as
\bnm
\sum_{m,n=0}^\infty \frac{(p,m+n)(2-p,m+n)(-1)^n x^{m+n}}{(p,m)(2-p,n)}
 L_{m}^{(p-1)}(-y)L_{n}^{(1-p)}(y)
=(1-4xy)^{-1/2};
\enm
In fact, Exton's result (3.13) can also be obtained from Exton's result (3.5)
by performing the replacements $d=1/2$ and $p'=2-p$.

Exton's result (4.3) should be read as
\bnm
\sum_{m,n=0}^\infty\frac{(p,m+n)(-1)^n x^{m+n}}{(p,m)(p,n)}
L_{m}^{(p-1)}(y)L_{n}^{(p-1)}(y)
={_0F_1}(-;p;-x^2y^2)
=\Gamma(p)(xy)^{1-p}J_{p-1}(2xy);
\enm

Exton's result (4.5) should be read as
\bnm
\sum_{m,n=0}^\infty\frac{(p,m+n)(2p-1,m+n)(-1)^n x^{m+n}}
{(p,m)(p,n)}L_{m}^{(p-1)}(y)L_{n}^{(p-1)}(y)
=(1+4x^2y^2)^{1/2-p};
\enm

Exton's result (5.3) should be read as
\bnm
\sum_{m,n=0}^\infty\frac{(1/2,m+n)(1/2,m+n)(-1)^m 2^{-2m-2n}x^{m+n}}
{(m+n)!(1/2,m)(1/2,n)m!\:n!}
H_{2m}(i y^{1/2})H_{2n}(y^{1/2})
=\text{exp}(4xy);
\enm

Exton's result (5.4) should be read as
\bnm
\sum_{m,n=0}^\infty\frac{(3/2,m+n)(2,m+n)(-1)^m 2^{-2-2m-2n}x^{m+n}}
{(m+n)!(3/2,m)(3/2,n)m!\:n!}
H_{2m+1}(i y^{1/2})H_{2n+1}(y^{1/2})
=i y \:\text{exp}(4xy);
\enm

Exton's result (5.5) should be read as
\bnm
\sum_{m,n=0}^\infty\frac{(3/2,m+n)(-1)^m 2^{-1-2m-2n}x^{m+n}}
{(1/2,m)(3/2,n)m!\:n!}
H_{2m}(i y^{1/2})H_{2n+1}(y^{1/2})
=y^{1/2}\text{exp}(4xy);
\enm

Exton's result (5.6) should be read as
\bnm
\sum_{m,n=0}^\infty\frac{(p',m+n)(p'-1/2,m+n)(-1)^{m+n} x^{m+n} 2^{-2m}}
{((2p'-1)/4,m+n)((2p'+1)/4,m+n)(1/2,m)(p',n)m!}\\
H_{2m}(y^{1/2})L_{n}^{(p'-1)}(-y)
=\text{cos}(4x^{1/2}y^{1/2});
\enm

Exton's result (5.7) should be read as
\bnm
\sum_{m,n=0}^\infty\frac{(1/2,m+n)(-1)^m x^{m+n}2^{-2m-2n}}
{(1/2,m)(1/2,n)m!\:n!}
H_{2m}(y^{1/2})H_{2n}(y^{1/2})
=\text{cos}(2xy);
\enm

Exton's result (5.8) should be read as
\bnm
\sum_{m,n=0}^\infty\frac{(3/2,m+n)(-1)^m x^{m+n+1}2^{-1-2m-2n}}
{(3/2,m)(3/2,n)m!\:n!}
H_{2m+1}(y^{1/2})H_{2n+1}(y^{1/2})
=\text{sin}(2xy);
\enm

Exton's result (6.2) should be read as
\bnm
\sum_{m=0}^q\frac{(-1)^m}
{(p,m)(p',q-m)}L_{m}^{(p-1)}(-y)L_{n}^{(p'-1)}(y) =\frac{((p+p'-1)/2,q)((p+p')/2,q)(-4y)^q}
{(p,q)(p',q)(p+p'-1,q)q!};
\enm

\section*{Acknowledgement}

X. Wang acknowledges support of National Natural Science Foundation of China (Grant No. 11201291),
and Natural Science Foundation of Shanghai, China (Grant No. 12ZR1443800).



\end{document}